\documentclass[10pt]{amsart}
\usepackage{amsmath, amsthm, amssymb}
\usepackage{mathbbol}
\usepackage{txfonts}

\newtheorem{thm}{Theorem}[section]
\newcommand{\bt}{\begin{thm}}
\newcommand{\et}{\end{thm}}

\newtheorem{cor}[thm]{Corollary}
\newcommand{\bc}{\begin{cor}}
\newcommand{\ec}{\end{cor}}

\newtheorem{lem}[thm]{Lemma}
\newcommand{\bl}{\begin{lem}}
\newcommand{\el}{\end{lem}}

\newtheorem{prop}[thm]{Proposition}
\newcommand{\bp}{\begin{prop}}
\newcommand{\ep}{\end{prop}}

\newtheorem{defn}[thm]{Definition}
\newcommand{\bd}{\begin{defn}}      
\newcommand{\ed}{\end{defn}}

\newtheorem{rmrk}[thm]{Remark}
\newcommand{\br}{\begin{rmrk}}
\newcommand{\er}{\end{rmrk}}

\newcommand{\thmref}[1]{Theorem~\ref{#1}}
\newcommand{\secref}[1]{Section~\ref{#1}}
\newcommand{\lemref}[1]{Lemma~\ref{#1}}

\newcommand{\corref}[1]{Corollary~\ref{#1}}
\newcommand{\propref}[1]{Proposition~\ref{#1}}

\newcommand{\N}{\mathbb{N}}

\newcommand{\R}{\mathbb{R}}
\newcommand{\Z}{\mathbb{Z}}

\newcommand{\dist}{\operatorname{dist}}
\newcommand{\diam}{\operatorname{diam}}

\newcommand{\hm}{{\mathcal H}}
\newcommand{\lm}{{\mathcal L}}
\newcommand{\md}{\operatorname{md}}

\newcommand{\jac}{{\mathbf J}}

\newcommand{\lip}{\operatorname{Lip}}
\newcommand{\mass}[2][]{{\mathbf M_{#1}}(#2)}

\newcommand{\form}{{\mathcal D}}        
\newcommand{\curr}{{\mathbf M}}         
\newcommand{\intrectcurr}{{\mathcal I}} 
\newcommand{\intcurr}{{\mathbf I}}      

\newcommand{\fillvol}{{\operatorname{Fillvol}}}
\newcommand{\fillrad}{{\operatorname{Fillrad}}}
\newcommand{\FV}{{\operatorname{FV}}}
\newcommand{\FR}{{\operatorname{FR}}}

\newcommand{\rstr}{\:\mbox{\rule{0.1ex}{1.2ex}\rule{1.1ex}{0.1ex}}\:}
\newcommand{\bdry}{\partial}

\newcommand{\on}[1]{|_{#1}}
\newcommand{\spt}{\operatorname{spt}}
\newcommand{\ohne}{\backslash}

\newcommand{\divergence}{{\operatorname{div}}}

\newcommand{\Arank}{{\operatorname{asrk}}}
\newcommand{\topdim}{{\operatorname{Topdim}}}

\begin{document}

\title[Isoperimetric inequalities and the asymptotic rank]{The asymptotic rank of metric spaces}

\author{Stefan Wenger}

\address
  {Department of Mathematics\\
University of Illinois at Chicago\\
851 S. Morgan Street\\
Chicago, IL 60607--7045}
\email{wenger@math.uic.edu}

\date{October 16, 2008}


\thanks{Partially supported by NSF grant DMS 0707009}



\begin{abstract}
 In this article we define and study a notion of asymptotic rank for metric spaces and show in our main theorem that for a large class of spaces, the asymptotic rank is characterized by the growth of the higher filling functions. For a proper, cocompact, simply-connected geodesic metric space of non-curvature in the sense of Alexandrov the asymptotic rank equals its Euclidean rank.
  \end{abstract}

\maketitle

\bigskip

\section{Introduction}\label{section:introduction}

For a proper cocompact ${\rm CAT}(0)$-space or non-positively curved Riemannian manifold $X$, the Euclidean rank is defined as the maximal $n\in\N$ for which Euclidean $n$-space $\R^n$ isometrically embeds into $X$. In this article we study the following generalization of the Euclidean rank. Let $(X, d)$ be an arbitrary metric space. A metric space $(Z, d_Z)$ is said to be an asymptotic subset of $X$ if there exist a sequence of subsets $Z_j\subset X$ and $r_j\to\infty$ such that $(Z_j, r_j^{-1}d)$ converges in the Gromov-Hausdorff sense to $(Z, d_Z)$.

\bd\label{definition:asymptotic-rank}
The asymptotic rank of $X$, denoted by $\Arank(X)$, is the supremum over $n\in\N$ for which 
there exists an asymptotic subset $Z$ of $X$ and a biLipschitz embedding $\varphi:K\to Z$ with $K\subset\R^n$ compact and $\lm^n(K)>0$.
\ed

Here $\lm^n$ denotes the Lebesgue measure on $\R^n$.
One may equivalently use asymptotic cones instead of asymptotic subsets. In \secref{section:asymptotic-rank}  it will be shown that $\Arank(X)$ is the supremum over $n\in\N$ for which there exists an $n$-dimensional normed space whose unit ball is an asymptotic subset of $X$. 
If $X$ is a Hadamard space, that is a complete metric space which is ${\rm CAT}(0)$, then  $\Arank(X)$ is the maximal geometric dimension of an asymptotic cone of $X$. If $X$ is moreover proper and cocompact then $\Arank(X)$
coincides with its Euclidean rank. This follows from work of Kleiner \cite{Kleiner-local-structure}. We refer to \secref{section:asymptotic-rank} below for details. 
For a general metric space $X$ we have the following properties:
 \begin{enumerate}
  \item $\Arank(X)\leq\sup\{\topdim(C): \text{$C\subset Z$ cpt, $Z$ an asymptotic subset of $X$}\}$;
  \item $\Arank(X))\geq \sup\{n\in\N: \text{$\exists \psi: \R^n\to X$ quasi-isometric}\}.$
\end{enumerate}
Clearly, if $X$ is a geodesic Gromov hyperbolic metric space then $\Arank(X)=1$ because every asymptotic subset of $X$ is isometric to a subset of a real tree. On the other hand, not all geodesic spaces with asymptotic rank $1$ are Gromov hyperbolic as the the countable wedge sum of  circles $S_j^1$ of length $j$ shows.
However, a geodesic metric space $X$ with $\Arank(X)=1$ which has a quadratic isoperimetric inequality for curves is Gromov hyperbolic, see \corref{corollary:gromov-hyp-arank=1}.

The main results of the present article characterize the asymptotic rank in terms of the growth of higher isoperimetric filling functions $\FV_{k+1}$ for a class of metric spaces including all Hadamard spaces. 
In the generality of (complete) metric spaces $X$, a suitable notion of $k$-chains and $k$-cycles is provided by the theory of integral currents, developed by Ambrosio and Kirchheim in \cite{Ambr-Kirch-curr}. Definitions will be given in \secref{section:preliminaries}. For $k\geq 0$, the space of integral $k$-currents in $X$ is denoted by $\intcurr_k(X)$. Given an element $T\in\intcurr_k(X)$, its mass is denoted by $\mass{T}$, its boundary (defined if $k\geq 1$) is an element of $\intcurr_{k-1}(X)$ and denoted by $\bdry T$. 
We recall that $X$ is said to admit an isoperimetric inequality of Euclidean type for $\intcurr_k(X)$ if there exists $D>0$ such that 
 \begin{equation}\label{equation:intro-isop-eucl}
  \fillvol(T)\leq D\mass{T}^{\frac{k+1}{k}}
 \end{equation}
 for every $T\in\intcurr_k(X)$ with $\bdry T=0$, where $\fillvol(T)$ is the least mass of an $S\in\intcurr_{k+1}(X)$ with $\bdry S=T$.  Similarly, $X$ is said to admit a cone type inequality for $\intcurr_k(X)$ if there exists $C>0$ such that 
 \begin{equation}\label{equation:intro-cone-ineq}
  \fillvol(T)\leq C\diam(\spt T)\mass{T}
 \end{equation}
 for every $T\in\intcurr_k(X)$ with $\bdry T=0$ and with bounded support $\spt T$. 
Important classes of spaces admitting cone type inequalities (in every dimension) include: Hadamard spaces, more generally, geodesic metric spaces with a convex metric, and metric spaces admitting a convex bicombing \cite{Wenger-GAFA}; furthermore, Riemannian manifolds without focal points. See \secref{section:diameter-vol-ineq-spaces} for more general classes and for details.

It was shown in \cite{Gromov-filling} and \cite{Wenger-GAFA} that metric spaces $X$ which admit cone type inequalities for $\intcurr_m(X)$, $m=1,\dots, k$, admit isoperimetric inequalities of Euclidean type.
The converse is true for $k=1$ for spaces $(X,d)$ in which, for some $Q>0$, any two points $x,x'$ in $X$ can be joined by a curve of length at most $Qd(x,x')$. For $k\geq 2$, it is not clear what the precise relationship between admitting isoperimetric inequalities of Euclidean type and admitting cone type inequalities for $\intcurr_m(X)$, $m=1,\dots,k$, is. Note that, for example,  \eqref{equation:intro-isop-eucl} is stronger than \eqref{equation:intro-cone-ineq} if $T$ has sufficiently small mass and large diameter.

We turn to the main results of the present paper.
For $k\in\N$ define the filling volume function $\FV_{k+1}$ on $X$ by
\begin{equation}\label{equation:fillvol-def}
 \FV_{k+1}(r):= \sup\{\fillvol(T): \text{$T\in\intcurr_k(X)$ with $\bdry T=0$ and $\mass{T}\leq r$}\}.
\end{equation}
Note that $X$ admits an isoperimetric inequality of Euclidean type for $\intcurr_k(X)$ if and only if $\FV_{k+1}(r)\leq Dr^{\frac{k+1}{k}}$
for all $r\geq 0$ and for some constant $D$.
Our main result can now be stated as follows.

\bt\label{Theorem:Subeuclidean-isoperimetric}
 Let $X$ be a complete metric space such that, for some $Q>0$, any two points $x,x'$ in $X$ can be joined by a curve of length at most $Qd(x,x')$. 
Let $k\in\N$ and suppose $X$ admits cone type inequalities for $\intcurr_m(X)$ for $m=1,\dots, k$. 
 If $k\geq \Arank(X)$ then
 \begin{equation}\label{equation:fillvol-function-sub-eucl-growth}
  \limsup_{r\to\infty}\frac{\FV_{k+1}(r)}{r^\frac{k+1}{k}}= 0.
 \end{equation}
In other words, $X$ admits a sub-Euclidean isoperimetric inequality for $\intcurr_k(X)$.
\et
\thmref{Theorem:Subeuclidean-isoperimetric} seems to be new even for cocompact Hadamard manifolds. For symmetric spaces of non-compact type a stronger result is known however: They admit linear isoperimetric inequalities above the rank, thus $\FV_{k+1}(r)\leq D r$ for all $r\geq 0$ for some constant $D$.
First applications of our theorem to the asymptotic geometry of non-positively curved spaces are given in \cite{Kleiner-Lang}, see also \cite{Kleiner-Lang-Oberwolfach}. 
A simple consequence of the above theorem is:

\bc\label{corollary:gromov-hyp-arank=1}
 Let $X$ be a geodesic metric space which admits a Euclidean isoperimetric inequality for $\intcurr_1(X)$, that is, has a quadratic isoperimetric inequality for curves. Then $X$ is Gromov hyperbolic if and only if $\Arank(X)=1$.
\ec

While one implication is clear, the other one follows from the theorem above and the well-known fact that geodesic metric spaces with a subquadratic isoperimetric inequality for curves are Gromov hyperbolic.

Our next theorem gives a lower bound for the growth of the filling functions below the asymptotic rank.

\bt\label{Theorem:lower-bound-fillvol}
 Let $X$ be a complete metric space such that, for some $Q>0$, any two points $x,x'$ in $X$ can be joined by a curve of length at most $Qd(x,x')$. Let $k\in\N$ and suppose $X$ admits isoperimetric inequalities of Euclidean type for $\intcurr_m(X)$ 
 with some constants $D_m$, $m=1,\dots, k-1$. 
 If $k< \Arank(X)$ then 
\begin{equation*}
 \FV_{k+1}(r)\geq\FV_{k+1}(X, L^\infty(X),r)\geq \varepsilon_k r^{\frac{k+1}{k}}
\end{equation*}
for all $r>0$ large enough. Here $\varepsilon_k>0$ is a suitable constant only depending on $D_m$, $m=1,\dots,k-1$. 
\et
Here $\FV_{k+1}(X,L^\infty(X),r)$ is defined analogously to $\FV_{k+1}(r)$ but the filling volume in $X$ is replaced by the filling volume in $L^\infty(X)$, see \secref{section:integral-currents}.
The proof of \thmref{Theorem:lower-bound-fillvol} will in fact show that not only the filling volume but also the filling radius function is bounded from below. 
As regards the constants $\varepsilon_k$ in the theorem, it can be shown that a geodesic metric space $X$ with $\Arank(X)>1$ satisfies
\begin{equation}\label{equation:opt-value-isop-const-curves}
 \limsup_{r\to\infty}\frac{\FV_2(X,L^\infty(X),r)}{r^2}\geq \frac{1}{4\pi}.
\end{equation}
In \cite{Wenger-Gromov-hyp-isop} it is moreover proved that a geodesic metric space $X$ which admits a (coarse) quadratic isoperimetric inequality for curves and for which 
\eqref{equation:opt-value-isop-const-curves} fails is Gromov hyperbolic and thus all its asymptotic subsets are isometric subsets of real trees.
As a consequence of Theorems \ref{Theorem:Subeuclidean-isoperimetric} and \ref{Theorem:lower-bound-fillvol} we have:

\bc
The higher isoperimetric filling functions $\FV_k$ detect the asymptotic rank of complete metric spaces which admit cone type inequalities and for which, for some $Q>0$, any two points can be joined by a curve of length at most $Q$ times their distance. In particular, they detect the Euclidean rank of cocompact Hadamard spaces.
\ec

Related results have been obtained for symmetric spaces of non-compact type (\cite{Brady-Farb}, \cite{Leuzinger}, \cite{Hindawi-divergence}) and for proper cocompact Hadamard spaces (\cite{Wenger-filling}) for the higher divergence invariants $\divergence_k$ of Brady and Farb.

As mentioned above, symmetric spaces $X$ of non-compact type admit linear isoperimetric inequalities for $\intcurr_k(X)$ for all $k\geq\Arank(X)$. It is an open problem, see \cite{Gromov-asymptotic}, whether this holds for more general spaces such as for example for all cocompact Hadamard manifolds. In \cite{Papasoglu-subeucl-linear} Papasoglu shows that for a simplicial complex $X$ with $H_1(X)=H_2(X)=0$ and for which every extremal $2$-cycle for $\FV_3$ has genus at most some fixed $g\in\N$, the following holds:
If $X$ admits a quadratic isoperimetric inequality for curves and a sub-Euclidean isoperimetric inequality for $2$-cycles then for every $\varepsilon>0$ there exists $D_\varepsilon$ such that
\begin{equation*}
 \FV_{3}(r)\leq D_\varepsilon r^{1+\varepsilon}
\end{equation*}
for every $r\geq 1$. In the notation of \cite{Gromov-metric-structures}, this says that $X$ admits an isoperimetric inequality of 
infinite (i.e.~arbitrary large) rank for $\intcurr_2(X)$. It seems to be unknown at present, whether in $\R^3$, endowed with a non-positively curved metric, extremal $2$-cycles have a uniform bound on their genus.

Our next theorem shows that, under suitable conditions, an isoperimetric inequality of infinite rank for $m$-cycles is passed on to $(m+1)$-cycles. More precisely:

\bt\label{theorem:infinite-rank-isop-progression}
 Let $X$ be a complete metric space, $k,k'\in\N$ with $k'\leq k$,
 and suppose that $X$ admits a cone type inequality for $\intcurr_m(X)$ for each
 $m=k',\dots,k$. If $X$ admits an isoperimetric inequality of infinite rank for $\intcurr_{k'}(X)$ then $X$ admits an isoperimetric inequality of
 infinite rank for $\intcurr_k(X)$.
\et

For example, if $X$ is geodesic and Gromov hyperbolic it can be shown that there exists a geodesic thickening $X_\varrho$ of $X$ which admits cone type 
inequalities for $\intcurr_m(X_\varrho)$ for all $m\geq 1$. In particular, $X_\varrho$ is Gromov hyperbolic and admits a linear isoperimetric inequality for 
$\intcurr_1(X_\varrho)$. Thus, by the theorem above, $X_\varrho$ admits an isoperimetric inequality of 
infinite rank for $\intcurr_k(X_\varrho)$ for every $k\geq 1$. The following problem seems to be unsolved:\\
\\
{\bf Problem:} Let $Y$ be geodesic and Gromov hyperbolic, $k\geq 2$, and suppose $Y$ admits cone type inequalities
 for $\intcurr_m(Y)$ for $m=1,\dots,k$. Is it true that $Y$ admits a linear isoperimetric inequality for $\intcurr_k(Y)$?\\
\\
This is known to be true under suitable conditions on the geometry on small scales, see \cite{Lang-hyperbolic}. The problem seems to be open however even in the case of Hadamard spaces that are Gromov hyperbolic but not proper and not cocompact.

\bigskip

One of the main ingredients in the proof of \thmref{Theorem:Subeuclidean-isoperimetric} is a `thick-thin' decomposition for integral currents which was proved in \cite{Wenger-compactness}, see also \secref{section:decomposition-theorem}.
This theorem can furthermore be used to establish polynomial isoperimetric inequalities on the large scale with exponent different from the Euclidean one for certain classes of metric spaces. We will do this at the end of this paper, in Sections \ref{section:proof-generalized-isop-inequaliy} and \ref{section:diameter-vol-ineq-spaces}.

\bigskip

{\bf Acknowledgments:} I would like to thank Bruce Kleiner, Urs Lang and Tim Riley for several discussions on topics related to this paper.

\section{Preliminaries}\label{section:preliminaries}

In this section we recall some definitions and facts that are used throughout the paper.

\subsection{Lipschitz maps and metric derivatives}
Let $(X,d_X)$ and $(Y, d_Y)$ be metric spaces. A map $\varphi: X\to Y$ is said to be Lipschitz continuous if there exists $C>0$ such that
\begin{equation*}
 d_Y(\varphi(x), \varphi(x'))\leq C d_X(x,x')\quad\text{for all $x,x'\in X$.}
\end{equation*}
 If there exists $C>0$ such that
\begin{equation*}
 C^{-1}d_X(x,x')\leq d_Y(\varphi(x), \varphi(x'))\leq C d_X(x,x')\quad\text{for all $x,x'\in X$}
\end{equation*}
 then $\varphi$ is said to be biLipschitz continuous.
Given a Lipschitz map $\varphi: U\to X$, where $U\subset\R^k$ is open, and given $v\in\R^k$, the metric directional derivative of $\varphi$ at $z\in U$ in direction $v$ is defined by
\begin{equation*}
 \md\varphi_z(v):= \lim_{r\searrow 0}\frac{d(\varphi(z+rv),\varphi(z))}{r}
\end{equation*}
if this limit exists. This notion of differentiability was first introduced and studied by Kirchheim in \cite{Kirchheim}, who in particular proved the following theorem. A similar statement was proved by Korevaar-Schoen in \cite{Korevaar-Schoen} around the same time. 
\bt
 Let $(X,d)$ be a metric space and $\varphi:U\to X$ a Lipschitz map, where $U\subset\R^n$ is open. Then for almost every $z\in U$ the metric directional derivative
 $\md\varphi_z(v)$ exists for every $v\in\R^n$. Furthermore, there are compact sets $K_i\subset U$, $i\in\N$, such that $\lm^n(U\ohne\cup K_i)=0$ and 
 such that the following property holds: 
 For every $i$ and every $\varepsilon>0$ there exists $r(i,\varepsilon)>0$ such that
 \begin{equation}\label{equation:strong-metric-norm}
  |d(\varphi(z+v),\varphi(z+w))-\md\varphi_z(v-w)|\leq \varepsilon|v-w|
 \end{equation}
 for all $z\in K_i$ and all $v,w\in\R^n$ satisfying $|v|,|w|\leq r(i,\varepsilon)$ and $z+w\in K_i$.
\et

Here $|\cdot|$ denotes the Euclidean norm.
If $\md\varphi_z(v)$ exists for all $v\in\R^n$ and satisfies \eqref{equation:strong-metric-norm} then $\md\varphi_z$  is called metric derivative of $\varphi$ 
at the point $z$. Clearly, the metric derivative is a seminorm, and a norm if $\varphi$ is biLipschitz.
It is not difficult to prove that if $U\subset\R^n$ is merely Borel measurable then $\md\varphi_z$ can be defined at almost every Lebesgue density point $z\in U$ 
by a simple approximation argument. The following is then an easy consequence of the above theorem and the remarks above.

\bc\label{corollary:metric-diff-bilipschitz}
 Let $Z$ be a metric space and $\varphi:K\to Z$ biLipschitz with $K\subset\R^n$ Borel measurable and such that $\lm^n(K)>0$. Then 
 there exists a norm $\|\cdot\|$ on $\R^n$ with the following property: For every $\varepsilon>0$ and for every finite set $S\subset\R^n$ there exist
 $r>0$ and a map $\psi: S\to Z$ such that $\psi: (S,r\|\cdot\|)\to Z$ is $(1+\varepsilon)$-biLipschitz.
\ec

This corollary will be used repeatedly in our paper.

\subsection{Integral currents in metric spaces}\label{section:integral-currents}
The theory of normal and integral currents in metric spaces was developed by Ambrosio and Kirchheim in \cite{Ambr-Kirch-curr} 
and provides a suitable framework for studying filling problems in the generality we wish to work in. 
In Euclidean space Ambrosio-Kirchheim's theory agrees with the classical theory of Federer-Fleming normal and integral currents \cite{Federer-Fleming}.

Let $(X,d)$ be a complete metric space and $k\geq 0$ and let $\form^k(X)$ be the set of $(k+1)$-tuples $(f,\pi_1,\dots,\pi_k)$ 
of Lipschitz functions on $X$ with $f$ bounded. The Lipschitz constant of a Lipschitz function $f$ on $X$ will
be denoted by $\lip(f)$.
\bd
A $k$-dimensional metric current  $T$ on $X$ is a multi-linear functional on $\form^k(X)$ satisfying the following
properties:
\begin{enumerate}
 \item If $\pi^j_i$ converges point-wise to $\pi_i$ as $j\to\infty$ and if $\sup_{i,j}\lip(\pi^j_i)<\infty$ then
       \begin{equation*}
         T(f,\pi^j_1,\dots,\pi^j_k) \longrightarrow T(f,\pi_1,\dots,\pi_k).
       \end{equation*}
 \item If $\{x\in X:f(x)\not=0\}$ is contained in the union $\bigcup_{i=1}^kB_i$ of Borel sets $B_i$ and if $\pi_i$ is constant 
       on $B_i$ then
       \begin{equation*}
         T(f,\pi_1,\dots,\pi_k)=0.
       \end{equation*}
 \item There exists a finite Borel measure $\mu$ on $X$ such that
       \begin{equation}\label{equation:mass-def}
        |T(f,\pi_1,\dots,\pi_k)|\leq \prod_{i=1}^k\lip(\pi_i)\int_X|f|d\mu
       \end{equation}
       for all $(f,\pi_1,\dots,\pi_k)\in\form^k(X)$.
\end{enumerate}
\ed
The space of $k$-dimensional metric currents on $X$ is denoted by $\curr_k(X)$ and the minimal Borel measure $\mu$
satisfying \eqref{equation:mass-def} is called mass of $T$ and written as $\|T\|$. We also call mass of $T$ the number $\|T\|(X)$ 
which we denote by $\mass{T}$.
The support of $T$ is, by definition, the closed set $\spt T$ of points $x\in X$ such that $\|T\|(B(x,r))>0$ for all $r>0$. 
Note that currents have by definition finite mass. Recently, a variant of Ambrosio-Kirchheim's theory that does not rely on the finite mass axiom (iii) has been developed by Lang in \cite{Lang-currents}. 

Every function $\theta\in L^1(K,\R)$ with $K\subset\R^k$ Borel measurable induces an element of $\curr_k(\R^k)$ by
\begin{equation*}
 \Lbrack\theta\Rbrack(f,\pi_1,\dots,\pi_k):=\int_K\theta f\det\left(\frac{\partial\pi_i}{\partial x_j}\right)\,d\lm^k
\end{equation*}
for all $(f,\pi_1,\dots,\pi_k)\in\form^k(\R^k)$.
The restriction of $T\in\curr_k(X)$ to a Borel set $A\subset X$ is given by 
\begin{equation*}
  (T\rstr A)(f,\pi_1,\dots,\pi_k):= T(f\chi_A,\pi_1,\dots,\pi_k).
\end{equation*}
This expression is well-defined since $T$ can be extended to a functional on tuples for which the first argument lies in 
$L^\infty(X,\|T\|)$.

If $k\geq 1$ and $T\in\curr_k(X)$ then the boundary of $T$ is the functional
\begin{equation*}
 \bdry T(f,\pi_1,\dots,\pi_{k-1}):= T(1,f,\pi_1,\dots,\pi_{k-1}).
\end{equation*}
It is clear that $\bdry T$ satisfies conditions (i) and (ii) in the above definition. If $\bdry T$ also satisfies (iii) then $T$ is called a normal current.
By convention, elements of $\curr_0(X)$ are also called normal currents.

The push-forward of $T\in\curr_k(X)$ 
under a Lipschitz map $\varphi$ from $X$ to another complete metric space $Y$ is given by
\begin{equation*}
 \varphi_\# T(g,\tau_1,\dots,\tau_k):= T(g\circ\varphi, \tau_1\circ\varphi,\dots,\tau_k\circ\varphi)
\end{equation*}
for $(g,\tau_1,\dots,\tau_k)\in\form^k(Y)$. This defines a $k$-dimensional current on $Y$.
It follows directly from the definitions that $\bdry(\varphi_{\#}T) = \varphi_{\#}(\bdry T)$.

\sloppy
We will mainly be concerned with integral currents. We recall that an $\hm^k$-measurable set $A\subset X$
is said to be countably $\hm^k$-rectifiable if there exist countably many Lipschitz maps $\varphi_i :B_i\longrightarrow X$ from subsets
$B_i\subset \R^k$ such that
\begin{equation*}
\hm^k\left(A\ohne \bigcup \varphi_i(B_i)\right)=0.
\end{equation*}
\fussy

An element $T\in\curr_0(X)$ is called integer rectifiable if there exist finitely many points $x_1,\dots,x_n\in X$ and $\theta_1,\dots,\theta_n\in\Z\ohne\{0\}$ such
that
\begin{equation*}
 T(f)=\sum_{i=1}^n\theta_if(x_i)
\end{equation*}
for all bounded Lipschitz functions $f$.
 A current $T\in\curr_k(X)$ with $k\geq 1$ is said to be integer rectifiable if the following properties hold:
 \begin{enumerate}
  \item $\|T\|$ is concentrated on a countably $\hm^k$-rectifiable set and vanishes on $\hm^k$-negligible Borel sets.
  \item For any Lipschitz map $\varphi:X\to\R^k$ and any open set $U\subset X$ there exists $\theta\in L^1(\R^k,\Z)$ such that 
    $\varphi_\#(T\rstr U)=\Lbrack\theta\Rbrack$.
 \end{enumerate}
Integer rectifiable normal currents are called integral currents. The corresponding space is denoted by $\intcurr_k(X)$. In case $X=\R^N$ is Euclidean space, $\intcurr_k(X)$ agrees with the space of $k$-dimensional Federer-Fleming integral currents in $\R^N$. If $A\subset\R^k$ is a Borel set of finite measure and
finite perimeter then $\Lbrack\chi_A\Rbrack \in\intcurr_k(\R^k)$. Here, $\chi_A$ denotes the characteristic function. If $T\in\intcurr_k(X)$ and if $\varphi:X\to Y$ is a Lipschitz 
map into another complete metric space then $\varphi_{\#}T\in\intcurr_k(Y)$. Moreover, every Lipschitz chain in a complete metric space $X$ can be viewed as an integral current in $X$.

\subsection{Generalized filling volume and filling radius functions}

Let $X,Y$ be complete metric spaces and suppose 
$X$ isometrically embeds into $Y$. 
Then the filling volume of $T\in\intcurr_k(X)$ in $Y$ is defined as
\begin{equation*}
 \fillvol_Y(T):= \inf\{\mass{S}: S\in\intcurr_{k+1}(X), \bdry S=T\}
\end{equation*}
where we agree on $\inf\emptyset=\infty$. In case $Y=X$ we have $\fillvol_Y(T)=\fillvol(T)$. Furthermore, for $r\geq 0$, we set
\begin{equation*}
 \FV_{k+1}(X,Y,r):= \sup\{\fillvol_Y(T): T\in\intcurr_k(X), \bdry T=0, \mass{T}\leq r\}.
\end{equation*}
Clearly, $\FV_{k+1}(X,X,r)=\FV_{k+1}(r)$ and furthermore
\begin{equation*}
 \fillvol_Y(T)\leq \fillvol(T)\quad\text{ and }\quad \FV_{k+1}(X,Y,r)\leq \FV_{k+1}(r).
\end{equation*}
The left hand sides of both inequalities are smallest for $Y:=L^\infty(X)$.
Here, $L^\infty(X)$ is the Banach space of bounded functions on $X$ with the supremum norm $$\|f\|_\infty:= \sup_{x\in X}|f(x)|.$$
Similarly, the filling radius of $T\in\intcurr_k(X)$ in $Y$ is defined as
\begin{equation*}
 \fillrad_Y(T):= \inf\{\varrho\geq 0: \text{$\exists S\in\intcurr_{k+1}(Y)$ with $\bdry S=T$, $\spt S\subset B(\spt T,\varrho)$}\}
\end{equation*}
and furthermore
\begin{equation*}
 \FR_{k+1}(X,Y,r):= \sup\{\fillrad_Y(T): T\in\intcurr_k(X), \bdry T=0, \mass{T}\leq r\}.
\end{equation*}
The same obvious inequalities as for the filling volume hold for the filling radius, namely
\begin{equation*}
 \fillrad_Y(T)\leq \fillrad(T)\quad\text{ and }\quad \FR_{k+1}(X,Y,r)\leq \FR_{k+1}(r),
\end{equation*}
and the left-hand sides of both inequalities are smallest for $Y:=L^\infty(X)$.

\section{Basic properties of the asymptotic rank of a metric space}\label{section:asymptotic-rank}

The following fact gives a somewhat more geometric definition of the asymptotic rank.

\bp\label{proposition:alt-def-asrk}
 Let $X$ be a metric space. Then
$\Arank(X)$ is the supremum over $n\in\N$ for which there exists an $n$-dim. normed space V, subsets $S_j\subset X$ and a 
sequence $R_j\to\infty$ such that $\frac{1}{R_j}S_j\to B(0,1)\subset V$ in the Gromov-Hausdorff sense, where $B(0,1)$ denotes the closed unit ball in $V$.
\ep

\begin{proof}
 This is an easy consequence of \corref{corollary:metric-diff-bilipschitz}.
\end{proof}

The properties of $\Arank(X)$ listed at the beginning of the article are a direct consequence of \propref{proposition:alt-def-asrk}. 
The main reason for using the terminology `asymptotic rank' is its relationship to the Euclidean rank in the case of proper cocompact Hadamard spaces. For an account on the theory of non-positively curved metric spaces and in particular Hadamard spaces, see e.g.\ \cite{Ballmann, Bridson-Haefliger, Burago-Burago-Ivanov}. The  
following result is a direct consequence of \propref{proposition:alt-def-asrk} above and of Theorems A, C and D  of Kleiner \cite{Kleiner-local-structure}.

\bt
 Let $X$ be a metric space.
 If $X$ is a Hadamard space then $\Arank(X)$ is the maximal geometric dimension of an asymptotic cone of $X$. 
 If $X$ is a proper cocompact length space with a convex metric then
 \begin{equation*}
  \Arank(X) = \sup\{n\in\N: \text{$\exists V$ $n$-dim. normed space and $\psi: V\to X$ isometric}\}.
 \end{equation*}
 In particular, if $X$ is a proper cocompact Hadamard space then $\Arank(X)$ equals its Euclidean rank.
\et

For the definition of geometric dimension see \cite{Kleiner-local-structure}.
Here a geodesic metric space $X$ is said to have a convex metric if for every pair of constant-speed geodesic segments $c_1,c_2:[0,1]\to X$ the function 
$t\mapsto d(c_1(t),c_2(t))$ is convex.

\section{A decomposition theorem for integral currents}\label{section:decomposition-theorem}
A crucial ingredient in the proof of \thmref{Theorem:Subeuclidean-isoperimetric} will be \thmref{theorem:suitable-decomposition} below, which gives a kind of `thick-thin' decomposition for integral currents. This theorem was proved in \cite{Wenger-compactness} and can in fact be used to furthermore establish polynomial isoperimetric inequalities and \thmref{theorem:infinite-rank-isop-progression} , see \secref{section:proof-generalized-isop-inequaliy}. We start with  the following definition.

\bd
 Let $k\geq 2$ and $\alpha>1$. A complete metric space $X$ is said to admit an isoperimetric inequality of rank $\alpha$ for $\intcurr_{k-1}(X)$ 
 if 
 there is a constant $D>0$ such that 
 \begin{equation}\label{equation:isop-ineq-alpha}
  \FV_k(r)\leq DI_{k,\alpha}(r)
 \end{equation}
 for all $r\geq 0$, 
 where $I_{k,\alpha}$ is the function given by
 \begin{equation*}
 I_{k,\alpha}(r):=\left\{\begin{array}{l@{\qquad}l}
  r^{\frac{k}{k-1}} & 0\leq r\leq 1\\
  r^{\frac{\alpha}{\alpha-1}} & 1 < r <\infty.
 \end{array}\right.
\end{equation*}
\ed
In \cite[6.32]{Gromov-metric-structures} the polynomial bound $r^{\frac{\alpha}{\alpha-1}}$ was termed an isoperimetric inequality 
{\it of rank greater than} $\alpha$.
Here we will use the shorter terminology {\it of rank} $\alpha$.
Isoperimetric inequalities of rank $k$ for $\intcurr_{k-1}(X)$ are exactly those of {\it Euclidean type}. 

Now, set $\Lambda:= \{(k,\alpha)\in\N\times(1,\infty): k\geq 2\}\cup\{(1,0)\}$, let $\gamma\in(0,\infty)$ and define auxiliary functions by
\begin{equation*}
 F_{1,0,\gamma}(r)=\gamma r\quad\text{and}\quad G_{1,0}(r)=r
\end{equation*}
and for $(k,\alpha)\in\Lambda\ohne\{(1,0)\}$
\begin{equation*}
 F_{k,\alpha,\gamma}(r):=\left\{\begin{array}{l@{\qquad}l}
  \gamma\cdot r^k & 0\leq r\leq 1\\
  \gamma\cdot r^\alpha & 1 < r <\infty
 \end{array}\right.
\end{equation*}
and
\begin{equation*}
 G_{k,\alpha}(r):=\left\{\begin{array}{l@{\qquad}l}
  r^{\frac{1}{k}} & 0\leq r\leq 1\\
  r^{\frac{1}{\alpha}} & 1 < r <\infty.
 \end{array}\right.
\end{equation*}

The thick-thin decomposition theorem alluded to above can now be stated as follows.

\bt\label{theorem:suitable-decomposition}
 Let $X$ be a complete metric space, $(k,\alpha)\in\Lambda$, and suppose in case $k\geq 2$ that $X$ admits an isoperimetric inequality of rank $\alpha$ for 
 $\intcurr_{k-1}(X)$. Then for every $\lambda\in(0, 1)$ there exists $\gamma\in(0,1)$ with the following property. 
 Abbreviate $F:= F_{k,\alpha,\gamma}$ and $G:= G_{k,\alpha}$ and let $\delta\in(0,1)$.
 For every $T\in\intcurr_k(X)$ there exist $R\in\intcurr_k(X)$ and $T_j\in\intcurr_k(X)$, $j\in\N$, such that
 \begin{equation*}
  T= R + \sum_{j=1}^\infty T_j
 \end{equation*}
 and for which the following properties hold:
 \begin{enumerate}
  \item $\bdry R=\bdry T$ and $\bdry T_j=0$ for all $j\in\N$;
  \item For all $x\in\spt R\ohne\spt\bdry T$ and all $0\leq r\leq \min\{5\delta G(\mass{R}), \dist(x,\spt\bdry T)\}$
   \begin{equation*}
    \|R\|(B(x,r))\geq \frac{1}{2}5^{-(k+\alpha)}F(r);
   \end{equation*}
  \item $\mass{T_j}\leq (1+\lambda)\nu\gamma\mass{T}$ for all $j\in\N$, where $\nu:=\delta$ if $k=1$ or $\nu:=\max\{\delta^k,\delta^\alpha\}$ otherwise;
  \item $\diam(\spt T_j)\leq 4G\left(\gamma^{-1}\frac{2}{1-\lambda}5^{k+\alpha}\mass{T_j}\right)$;
  \item $\mass{R}+\frac{1-\lambda}{1+\lambda}\sum_{i=1}^\infty\mass{T_i}\leq \mass{T}$.
 \end{enumerate} 
\et
If $k=1$, all statements of the theorem hold for $\lambda=0$ as well. For the proof of \thmref{Theorem:Subeuclidean-isoperimetric} we will only need the case $\alpha = k$. The general case will be used to prove new polynomial isoperimetric inequalities in \secref{section:proof-generalized-isop-inequaliy}. For the proof of \thmref{theorem:suitable-decomposition} we refer to \cite{Wenger-compactness}.

We will furthermore need the following proposition from \cite{Wenger-compactness}, see also \cite{Ambr-Kirch-curr}.

\bp\label{proposition:isop-ineq-growth-estimate}
 Let $X$ be a complete metric space, $k\geq 2$, $\alpha>1$, and suppose that $X$ admits an isoperimetric inequality of rank 
 $\alpha$ for $\intcurr_{k-1}(X)$ with a constant $D_{k-1}\in [1,\infty)$.
 Then for every $T\in\intcurr_{k-1}(X)$ with $\bdry T=0$ and every $\varepsilon>0$ there exists an $S\in\intcurr_k(X)$ with $\bdry S=T$, satisfying
 \begin{equation}\label{equation:upper-bound-mass-volgrowth-estimate-prop}
  \mass{S}\leq \min\left\{(1+\varepsilon)\fillvol(T), D_{k-1}I_{k,\alpha}(\mass{T})\right\}
 \end{equation}
 and with the following property: For every $x\in \spt S$ and every $0\leq r\leq \dist(x, \spt T)$ we have
  \begin{equation*}
  \|S\|(B(x,r))\geq F_{k,\alpha,\mu}(r)
 \end{equation*}
 where 
 \begin{equation*}
  \mu:= \min\left\{\frac{1}{(3D_{k-1})^{k-1}\alpha_1^k}, \frac{1}{(3D_{k-1})^{\alpha-1}\alpha_1^\alpha}\right\}
 \end{equation*}
 with $\alpha_1:=\max\{k,\alpha\}$.
\ep

A direct consequence of the proposition is the following estimate on the filling radius.

\bc\label{corollary:fillrad-fillvol}
 Let $X$ be a complete metric space, $k\geq 2$, $\alpha>1$, and suppose that $X$ admits an isoperimetric inequality of rank $\alpha$ for $\intcurr_{k-1}(X)$. 
 Then for every $T\in\intcurr_{k-1}(X)$ with $\bdry T=0$ we have
 \begin{equation*}
  \fillrad_X(T)\leq G_{k,\alpha}(\mu^{-1}\fillvol_X(T))\leq\left\{\begin{array}{l@{\qquad}l}
   \mu'\mass{T}^{\frac{1}{k-1}} &\mass{T}\leq 1\\
   \mu'\mass{T}^{\frac{1}{\alpha-1}} &\mass{T}>1,
  \end{array}\right.
 \end{equation*}
 where 
 \begin{equation*}
  \mu':=\max\left\{\left(\frac{D_k}{\mu}\right)^{\frac{1}{k}}, \left(\frac{D_k}{\mu}\right)^{\frac{1}{\alpha}}\right\}.
 \end{equation*}
\ec

We end this section with the following useful fact.

\bl\label{lemma:optimal-decomposition-R}
 Let $X$ be a complete metric space, $(k,\alpha)\in\Lambda$, and $0<\varepsilon,\delta\leq 1$. Set $F:= F_{k,\alpha, \varepsilon}$ and $G:= G_{k,\alpha}$ 
 and let $R\in\intcurr_k(X)$ be such that $\bdry R = 0$ and
 \begin{equation*}
  \|R\|(B(x,r))\geq F(r)
 \end{equation*}
 for all $x\in\spt R$ and all $r\in[0,\delta G(\mass{R})]$.
 Then there exist constants $m\in\N$ and $E>0$ depending only on $k$, $\alpha$, $\delta$, $\varepsilon$
 and a decomposition $R=R_1+\dots+R_m$ with $R_i\in\intcurr_k(X)$, $\bdry R_i=0$, and 
 \begin{enumerate}
  \item $\|R_i\|(B(x,r))\geq F(r)$ for all $x\in\spt R_i$ and all $r\in[0,\delta G(\mass{R})]$
  \item $\mass{R} = \mass{R_1}+\dots+\mass{R_m}$
  \item $\diam(\spt R_i)\leq E G(\mass{R_i})$.
 \end{enumerate}
\el

It follows from (i) that, in particular, 
\begin{equation*}
 \mass{R_i}\geq \varepsilon\min\{\delta^k, \delta^\alpha\}\, \mass{R}.
\end{equation*}
\begin{proof}
 Set $\alpha':= \alpha$ if $k\geq 2$ or $\alpha':= 1$ if $k=1$.
 Fix $x\in\spt R$ arbitrary and observe that
\begin{equation*}
 \|R\|\left(B\left(x,t+2^{-1}\delta G(\mass{R})\right)\Big\backslash B\left(x,t-2^{-1}\delta G(\mass{R})\right)\right)=0
\end{equation*}
 for some $t\in \left[\frac{5}{2}\delta G(\mass{R}),\frac{3}{\varepsilon}\max\{\delta^{1-k},\delta^{1-\alpha'}\}G(\mass{R})\right]$.
Thus $R_1:= R\rstr B(x,t)$ satisfies $R_1\in\intcurr_k(X)$, $\bdry R_1=0$ and
\begin{equation}\label{equation:growth-equi-compactness-boundary}
 \|R_1\|(B(x',r))\geq F(r)
\end{equation}
for all $x'\in\spt R_1$ and all $0\leq r\leq \delta G(\mass{R})$.
In particular, we have
\begin{equation}\label{equation:mass-bound-opt-decomp}
 \mass{R_1}\geq \varepsilon \min\{\delta^k, \delta^{\alpha'}\} \mass{R}
\end{equation}
and thus
\begin{equation*}
 \diam(\spt R_1)\leq EG(\mass{R_1})
\end{equation*}
for a constant $E$ depending only on $k, \alpha,\delta,\varepsilon$.
Proceeding in the same way with $R- R_1$ one eventually obtains a decomposition $R=R_1+\dots+R_m$ with
the desired properties. The bound on $m$ clearly follows from \eqref{equation:mass-bound-opt-decomp}.
\end{proof}

\section{Sub-Euclidean isoperimetric inequalities}\label{section:proof-subeuclidean-isop}

In this section we prove the main result of this paper, \thmref{Theorem:Subeuclidean-isoperimetric}. We begin with the following simple lemma.
\bl\label{lemma:power-sum}
 Let $k\geq 2$, $\alpha>1$ and $0<\lambda, \delta \leq 1$. If $L>0$ and $0\leq t_i < \delta L$ are such that
 \begin{equation*}
  \lambda\sum_{i=1}^\infty t_i\leq L
 \end{equation*}
 then
 \begin{equation*}
  \sum_{i=1}^\infty I_{k,\alpha}(t_i) 
   \leq \frac{2(1+\delta\lambda)}{\lambda}\max\left\{(2\delta)^{\frac{1}{k-1}},(2\delta)^{\frac{1}{\alpha-1}}\right\}I_{k,\alpha}(L).
 \end{equation*}
\el
\begin{proof}
 Pick finitely many integer numbers $0=:m_0<m_1<m_2<\dots<m_{j_0}$ with the property that
 \begin{equation*}
  \delta L< t_{m_{i-1}+1}+\dots+t_{m_i}<2\delta L
 \end{equation*}
 for each $i=1,\dots,j_0$ and
 \begin{equation*}
  \sum_{n=m_{j_0}+1}^\infty t_n\leq \delta L.
 \end{equation*}
 Then $j_0\leq \frac{1}{\lambda\delta}$ and hence
 \begin{equation*}
  \begin{split}
   \sum_{i=1}^\infty I_{k,\alpha}(t_i) &\leq \sum_{i=1}^{j_0}I_{k,\alpha}(t_{m_{i-1}+1}+\dots+t_{m_i}) + I_{k,\alpha}\left(\sum_{n=m_{j_0}+1}^\infty t_n\right)\\
        &\leq \frac{1}{\lambda\delta}I_{k,\alpha}(2\delta L) + I_{k,\alpha}(\delta L)\\
        &\leq \frac{2(1+\delta\lambda)}{\lambda}\max\left\{(2\delta)^{\frac{1}{k-1}},(2\delta)^{\frac{1}{\alpha-1}}\right\} I_{k,\alpha}(L).
  \end{split}
 \end{equation*}
\end{proof}

\bl\label{lemma:good-cycle-fillvol}
 Let $X$ be a complete metric space, $k\geq 1$, $\alpha>1$, and suppose $X$ admits an isoperimetric inequality of rank $\alpha$ for $I_k(X)$. 
 In case $k\geq 2$ suppose furthermore that $X$ also admits an isoperimetric inequality of rank $k$ for $\intcurr_{k-1}(X)$.
 Let $\varepsilon>0$ and $T\in\intcurr_k(X)$ with $\bdry T=0$. If $\fillvol(T)\geq\varepsilon I_{k+1,\alpha}(\mass{T})$ then there exists $T'\in\intcurr_k(X)$ 
 with $\bdry T'=0$ and satisfying the following properties:
 \begin{enumerate}
  \item $\fillvol(T')\geq \frac{\varepsilon}{2} I_{k+1,\alpha}(\mass{T'})$
  \item $\mass{T'}\geq A\mass{T}$
  \item $\diam(\spt T')\leq B\mass{T'}^{\frac{1}{k}}$
  \item $\|T'\|(B(x,r))\geq Cr^k$ for all $r\in[0,5\delta\mass{T'}^{\frac{1}{k}}]$.
 \end{enumerate}
 Here, $A,B,C,\delta>0$ are constants depending only on $k,\alpha,\varepsilon$ and the constants of the isoperimetric inequalities. 
\el

\begin{proof}
 Set $\lambda:=1/3$ and
 \begin{equation*}
  \delta:= \min\left\{\frac{3}{8},\frac{\varepsilon}{64D_k}, \left(\frac{\varepsilon}{64D_k}\right)^{\frac{\alpha -1}{k}}\right\},
 \end{equation*}
 where $D_k$ is the constant for the isoperimetric  inequality for $\intcurr_k(X)$. 
 Let $T=R+\sum_{j=1}^\infty T_j$ be a decomposition as in \thmref{theorem:suitable-decomposition}.
 It then follows from \lemref{lemma:power-sum} that
 \begin{equation*}
  \begin{split}
   \fillvol(T)&\leq \fillvol(R)+ D_k\sum_{i=1}^\infty I_{k+1,\alpha}(\mass{T_j})\\
              &\leq \fillvol(R)+ \frac{\varepsilon}{2}I_{k+1,\alpha}(\mass{T})
  \end{split}
 \end{equation*}
 and thus
 \begin{equation}\label{equation:suitable-fillvol-R}
  \fillvol(R)\geq \frac{\varepsilon}{2}I_{k+1,\alpha}(\mass{T})\geq \frac{\varepsilon}{2} I_{k+1,\alpha}(\mass{R}).
 \end{equation}
 This together with the isoperimetric inequality for $\intcurr_k(X)$ yields
 \begin{equation*}
  \mass{R}
   \geq \min\left\{\left(\frac{\varepsilon}{2D_k}\right)^{\frac{k}{k+1}},\left(\frac{\varepsilon}{2D_k}\right)^{\frac{\alpha-1}{\alpha}}\right\} \mass{T}.
 \end{equation*}
 Let $R=R_1+\dots+R_m$ be a decomposition of $R$ as in \lemref{lemma:optimal-decomposition-R}. By \eqref{equation:suitable-fillvol-R} and the special properties
 of the decomposition there exists an $i$ such that $T':= R_i$ satisfies
 \begin{equation*}
  \fillvol(T')\geq \frac{\varepsilon}{2}I_{k+1,\alpha}(\mass{T'}).
 \end{equation*}
 It is clear from \lemref{lemma:optimal-decomposition-R} that $T'$ satisfies all the desired properties.
\end{proof}

We are now ready for the proof of the sub-Euclidean isoperimetric inequality.

\begin{proof}[{Proof of \thmref{Theorem:Subeuclidean-isoperimetric}}]
We argue by contradiction and suppose therefore that
\begin{equation*}
 \limsup_{r\to\infty} \frac{\FV_{k+1}(r)}{r^{\frac{k+1}{k}}}\geq 2\varepsilon_0>0
\end{equation*}
for some $\varepsilon_0>0$.
In particular, there is a sequence $T_m\in\intcurr_k(X)$ with $\bdry T_m=0$ and such that $\mass{T_m}\to\infty$ and
\begin{equation}\label{equation:euclidean-behavior-fillvol}
 \fillvol(T_m)\geq \varepsilon_0\mass{T_m}^{\frac{k+1}{k}}
\end{equation}
for every $m\in\N$.
By Theorem 1.2 of \cite{Wenger-GAFA}, $X$ admits an isoperimetric inequality of Euclidean type for $\intcurr_k(X)$ and, if $k\geq 2$ also one for 
$\intcurr_{k-1}(X)$. Therefore we may assume by \lemref{lemma:good-cycle-fillvol} that
\begin{equation}\label{equation:sequence-diam}
 \diam(\spt T_m)\leq B\mass{T_m}^{\frac{1}{k}}
\end{equation}
and
\begin{equation}\label{equation:sequence-growth}
 \|T_m\|(B(x,r))\geq Cr^k
\end{equation}
for all $x\in\spt T_m$ and all $r\in[0,5\delta\mass{T_m}^{1/k}]$, where $B,C,\delta$ are constants independent of $m$.
We set $r_m:= \mass{T_m}^\frac{1}{k}$ and note that $r_m\to\infty$. We choose $S_m\in\intcurr_{k+1}(X)$ with $\bdry S_m=T_m$ and
\begin{equation*}
 \mass{S_m}\leq D_k[\mass{T_m}]^{\frac{k+1}{k}}
\end{equation*}
and with the volume growth property of \propref{proposition:isop-ineq-growth-estimate}. 
We define a sequence of metric spaces $X_m:=(X,\frac{1}{r_m}d_X)$ where $d_X$ denotes the metric on $X$. 
Setting $Z_m:= \spt S_m\subset X_m$ it follows directly from \propref{proposition:isop-ineq-growth-estimate} and \eqref{equation:sequence-diam} and
\eqref{equation:sequence-growth} that the sequence $(Z_m,\frac{1}{r_m}d_X)$ is uniformly compact.
Therefore, by Gromov's compactness theorem \cite{Gromov-poly-growth} there exists (after passage to a subsequence) a compact metric space $(Z,d_Z)$ and isometric embeddings 
$\varphi_m: (Z_m, \frac{1}{r_m}d_X)\hookrightarrow (Z,d_Z)$ and such that $\varphi_m(Z_m)$ is a Cauchy sequence with respect to the Hausdorff distance. 
Denote by $S'_m$ the current $S_m$ viewed as an element of $\intcurr_{k+1}(X_m)$. 
Since $\mass{\varphi_{m\#}S'_m}\leq D_k$ and $\mass{\bdry(\varphi_{m\#}S'_m)}=1$ we may assume by the compactness and closure theorems for currents 
that $\varphi_{m\#}S'_m$ weakly converges to some $S\in\intcurr_{k+1}(Z)$. 
We first show that $\bdry S\not=0$. For this we choose $x_m\in\spt S'_m$ arbitrarily and define an auxiliary metric space $Y$ as the disjoint union 
$\bigsqcup_{m=1}^\infty X_m$ and endow it with the metric $d_Y$ in such a way that $d_Y\on{X_m\times X_m}=\frac{1}{r_m}d_X$ as well as 
\begin{equation*}
d_Y(y,y')=\frac{1}{r_m}d_X(y, x_m)+3+\frac{1}{r_{m'}}d_X(y',x_{m'}) 
\end{equation*}
whenever $y\in X_m$ and $y'\in X_{m'}$ with $m'\not=m$. It is clear that $Y$ admits a local cone type inequality for 
$\intcurr_l(Y)$, $l=1,\dots, k$, in the sense of \cite{Wenger-flatconv} and that any two points in $Y$ at distance less than $2$ can be joined by a curve of length at most $Q$ times their distance. 
Denote by $T'_m$ the current $T_m$ viewed as an element of $\intcurr_k(Y)$ and note that $\mass{T'_m}=1$.
Now $T'_m$ cannot weakly converge to $0$ since otherwise, by Theorem 1.4 in \cite{Wenger-flatconv}, we have
$\fillvol(T'_m)\to 0$ and, in particular, there exist $\hat{S}_m\in\intcurr_{k+1}(Y)$ with $\bdry\hat{S}_m=T'_m$ for all $m\in\N$ and
such that $\mass{\hat{S}_m}\to 0$. Of course, it is not restrictive to assume that $\spt\hat{S}_m\subset X_m$. 
Denote by $\tilde{S}_m$ the current $\hat{S}_m$ viewed as a current in $X$. Then $\tilde{S}_m$ satisfies 
$\bdry\tilde{S}_m=T_m$ and
\begin{equation*}
 \frac{\mass{\tilde{S}_m}}{r_m^{k+1}}=\mass{\hat{S}_m}\to 0,
\end{equation*}
which contradicts \eqref{equation:euclidean-behavior-fillvol}.
Thus, $T'_m$ does not weakly converge to $0$ and therefore there exist $\varepsilon>0$ and Lipschitz maps 
$f,\pi_1,\dots,\pi_k\in\lip(Y)$ with $f$ bounded such that 
\begin{equation*}
 T'_m(f,\pi_1,\dots,\pi_k)\geq \varepsilon\quad\text{for all $m\in\N$.}
\end{equation*}
Note that $\cup Z_m\subset Y$ is bounded so that the functions $\pi_i$ are bounded on $\cup Z_m$. 
We define Lipschitz functions $f_m$ and $\pi_i^m$ on $\varphi_m(Z_m)$ by
$f_m(z):= f(\varphi_m^{-1}(z))$ and $\pi_i^m(z):=\pi_i(\varphi_m^{-1}(z))$ for $z\in\varphi_m(Z_m)$. Here, we view $\varphi^{-1}_m$ as a map from $\varphi(Z_m)$ 
to $Y=\sqcup_{l=1}^\infty X_l$ with image in $X_m\subset Y$. 
By McShane's extension theorem there exist extensions $\hat{f}_m, \hat{\pi}_i^m: Z\to\R$ of $f_m$ and $\pi_i^m$ with the same Lipschitz constants
as $f$ and $\pi_i$.
By Arzel\`a-Ascoli theorem we may assume that $\hat{f}_m$ and 
$\hat{\pi}_i^m$ converge uniformly to Lipschitz maps $\hat{f}$, $\hat{\pi}_i$ on $Z$. Finally, we abbreviate $T''_m:=\varphi_{m\#}T'_m$ and use 
\cite[Proposition 5.1]{Ambr-Kirch-curr} to estimate
\begin{align*}
   \bdry S (\hat{f},\hat{\pi}_1,\dots,\hat{\pi}_k) &= \lim_{m\to\infty}T''_m(\hat{f},\hat{\pi}_1,\dots,\hat{\pi}_k)\\
&=\lim_{m\to\infty}\Big[T''_m(\hat{f}_m,\hat{\pi}_1^m,\dots,\hat{\pi}_k^m) + T''_m(\hat{f}-\hat{f}_m,\hat{\pi}_1,\dots,\hat{\pi}_k)\\
  &\quad\qquad\qquad+T''_m(\hat{f}_m,\hat{\pi}_1,\dots,\hat{\pi}_k)-T''_m(\hat{f}_m,\hat{\pi}_1^m,\dots,\hat{\pi}_k^m)\Big]\\
  &\geq \varepsilon - \limsup_{m\to\infty}\Big[\prod_{i=1}^k\lip(\hat{\pi}_i)\int_Z|\hat{f}-\hat{f}_m|\,d\|T''_m\|\Big]\\
  & \qquad- \limsup_{m\to\infty}\Big[\lip(\hat{f}_m)\sum_{i=1}^k\int_Z|\hat{\pi}_i-\hat{\pi}_i^m|\,d\|T''_m\|\Big]\\
  &= \varepsilon.
\end{align*}
This shows that indeed $\bdry S\not=0$ and hence also $S\not=0$. 
Now, since the Hausdorff limit $Z':=\lim_H\varphi_m(Z_m)\subset Z$ is an asymptotic subset of $X$ and since $\spt S\subset Z'$ and $S\not=0$, Theorem 4.5 in \cite{Ambr-Kirch-curr} shows that there exists a biLipschitz map $\nu:K\subset\R^{k+1}\to Z'$ where
$K$ is compact and of strictly positive Lebesgue measure.
This is in contradiction with the hypothesis that $k\geq\Arank X$ and hence this completes the proof.
\end{proof}

The arguments in the proof above can easily be used to establish the following result.

\bt\label{Theorem:fillvol-Linfty-above-rank}
 Let $k\in\N$ and let $X$ be a complete metric space which admits an isoperimetric inequality of Euclidean type for $\intcurr_k(X)$ and, in case that  $k\geq 2$, also one for $\intcurr_{k-1}(X)$. 
 If $k\geq \Arank(X)$ then
 \begin{equation*}
  \limsup_{r\to\infty}\frac{\FV_{k+1}(X, L^\infty(X), r)}{r^\frac{k+1}{k}}= 0.
 \end{equation*}
\et

Note that in contrast to the main theorem we do not assume here that $X$ admits cone type inequalities.

\section{Lower bounds on the filling radius}\label{section:proof-lower-bound-fillrad}
In this section we prove the following theorem.

\bt\label{Theorem:lower-bound-fillrad}
 Let $X$ be a complete metric space such that, for some $D_0>0$, any two points $x,x'$ in $X$ can be joined by a curve of length at most $D_0d(x,x')$. Let $k\in\N$ and suppose $X$ admits isoperimetric inequalities of Euclidean type for $\intcurr_m(X)$ 
 with some constants $D_m$, $m=1,\dots, k-1$. 
 If $k< \Arank(X)$ then 
\begin{equation*}
 \FR_{k+1}(r)\geq\FR_{k+1}(X, L^\infty(X),r)\geq \varepsilon_k r^{\frac{1}{k}}
\end{equation*}
for all $r>0$ large enough and for some $\varepsilon_k>0$ depending only on $D_m$, $m=1,\dots,k-1$. 
\et

Note that \thmref{Theorem:lower-bound-fillvol} is a consequence of the above since, by \corref{corollary:fillrad-fillvol}, we have
\begin{equation*}
 \fillrad_{L^\infty(X)}(T)\leq C[\fillvol_{L^\infty(X)}(T)]^{\frac{1}{k+1}}
\end{equation*}
for some constant $C$ and all $T\in\intcurr_k(X)$ with $\bdry T=0$.

\begin{proof}
 Let $Z$ be an  asymptotic subset of $X$ and $\varphi: K\subset\R^{k+1}\to Z$ a biLipschitz map with $K$ compact and 
 such that $\lm^{k+1}(K)>0$.
 Let $\|\cdot\|$ be a norm on $\R^{k+1}$ as in \corref{corollary:metric-diff-bilipschitz} and set $V:=(\R^{k+1}, \|\cdot\|)$. Let $\{v_1,\dots,v_{k+1}\}\subset V$ 
 and  $\{v^*_1,\dots,v^*_{k+1}\}\subset V^*$ be bases satisfying
 \begin{equation*}
  \|v_i\| = 1 = \|v^*_i\|\quad\text{ and }\quad v^*_i(v_j) = \delta_{ij}\quad \text{for all {i,j}.}
 \end{equation*}
  
 Let $Q$ denote the cube $Q:= \left\{\sum_{i=1}^{k+1}\lambda_iv_i:0\leq \lambda_i\leq 1\right\}$.
For $n=1,\dots, k+1$ denote by $A(n)$ the set of increasing functions
\begin{equation*}
 \alpha: \{1,\dots,n\}\to\{1,\dots,k+1\},
\end{equation*}
and by $L_\alpha$ the subspace generated by $\{v_{\alpha(1)},\dots,v_{\alpha(n)}\}$ whenever $\alpha\in A(n)$. 
It is clear that 
\begin{equation*}
 1\leq \mu_{L_\alpha}^{m*}(v_{\alpha(1)}\wedge\dots\wedge v_{\alpha(n)})\leq n^{\frac{n}{2}},
\end{equation*}
where $\mu_{L_\alpha}^{m*}$ denotes the Gromov mass$^*$-volume on $L_\alpha$ and 
 $v_{\alpha(1)}\wedge\dots\wedge v_{\alpha(n)}$ is the
 parallelepiped spanned by these vectors. Recall that for a compact set $K\subset W$ in an $l$-dimensional normed space the associated integer rectifiable 
 current $\Lbrack K\Rbrack\in\intrectcurr_l(W)$ satisfies $\|\Lbrack K\Rbrack\|(A)=\mu_W^{m*}(A)$, see e.g.~\cite[Proposition 2.7]{Wenger-Gromov-hyp-isop}.
Fix $m\in\N$ large enough (as chosen below) and let $\mathcal{Q}_n$ denote the $n$-skeleton of the cubical subdivision of $Q$ given by
\begin{equation*}
 \mathcal{Q}_0:= Q\cap \left\{2^{-m}\sum_{i=1}^{k+1}\delta_iv_i: \delta_i\in\Z\right\}
\end{equation*}
if $n=0$ and
\begin{equation*}
 \mathcal{Q}_n:= \left\{\sigma= z+ 2^{-m}(L_\alpha\cap Q): z\in\mathcal{Q}_0, \alpha\in A(n), \sigma\subset Q\right\}
\end{equation*}
if $n\in\{1,\dots, k+1\}$. We furthermore set $\bdry \mathcal{Q}_n:= \{\sigma\in \mathcal{Q}_n: \sigma\subset\bdry Q\}$.
Let $s>0$ and $\varepsilon>0$ be arbitrary. By \corref{corollary:metric-diff-bilipschitz} there exists an $r>0$ and a $(1+\varepsilon)$-biLipschitz map
 $\hat{\psi}:(\bdry\mathcal{Q}_0, r\|\cdot\|)\to Z.$
The definition of asymptotic subset and the fact that $\bdry\mathcal{Q}_0$ is a finite set imply the existence of $s'\geq \max\{s, 2^{m+2}\}$ and a 
$(1+2\varepsilon)$-biLipschitz map $\psi: (\bdry\mathcal{Q}_0, s'\|\cdot\|)\to X$. 
We write $V':= (\R^{k+1}, s'\|\cdot\|)$ and note that by McShane's extension theorem, there exists a $(1+2\varepsilon)(k+1)$-Lipschitz extension 
$\eta: L^\infty(X)\to V'$ of the map $\psi^{-1}: \psi(\bdry\mathcal{Q}_0)\to V'$. In the following, we regard $Q, \mathcal{Q}_n$, 
and $\bdry\mathcal{Q}_n$ as subsets of $V'$.
We associate with $\bdry\mathcal{Q}_n$ the additive subgroup
\begin{equation*}
 G_n = \left\{T\in\intcurr_n(V'): T=\sum c_i\Lbrack\sigma_i\Rbrack,c_i\in\Z, \sigma_i\in\bdry\mathcal{Q}_n\right\}\subset\intcurr_n(V').
\end{equation*}
%
%
We now construct homomorphisms $\Lambda_n: G_n\to \intcurr_n(X)$ and  $\Gamma_n:G_n\to \intcurr_{n+1}(V')$ for $n= 0,1,\dots, k$, with the property that
for all $T\in G_n$
\begin{enumerate}
 \item $\bdry\circ\Lambda_n= \Lambda_{n-1}\circ\bdry$ whenever $n\geq 1$
 \item $\mass{\Lambda_n(T)}\leq C'_n\mass{T}$
 \item $\bdry \Gamma_n(T) = T-\eta_\#\Lambda_n(T) - \Gamma_{n-1}(\bdry T)$ whenever $n\geq 1$
 \item $\mass{\Gamma_n(T)}\leq 2^{-m}D'_n s'\mass{T}$
 \item $\spt(\eta_\#\Lambda_n(T))\subset B(\spt T, 2^{-m}E_n s')$
 \item $\spt(\Gamma_n(T))\subset B(\spt T, 2^{-m}E'_n s')$.
\end{enumerate}
Here, $C'_n,D'_n,E_n,E'_n$ are constants depending only on $D_0,\dots, D_n$.
For $n=0$ we simply set $\Lambda_0(\Lbrack x\Rbrack):= \psi_\#\Lbrack x\Rbrack$ for each $x\in\bdry\mathcal{Q}_0$ and extend $\Lambda_0$ to $G_0$ as a 
homomorphism. Set furthermore $\Gamma_0:= 0$ and note that for $n=0$ the above properties are satisfied with $C'_0:=1, D'_0:=0, E_0:=0, E'_0:=0$.
Suppose now that $\Lambda_{n-1}$ and $\Gamma_{n-1}$ have been defined for some $n\in\{1,\dots, k\}$ and that they have the properties listed above.
In order to define $\Lambda_n$ let $\sigma\in\bdry\mathcal{Q}_n$ and note that 
\begin{equation*}
 (2^{-m}s')^n \leq \mass{\Lbrack\sigma\Rbrack}\leq n^{\frac{n}{2}}(2^{-m}s')^n. 
\end{equation*}
If $n=1$ there exists 
a Lipschitz curve of length at most $(1+2\varepsilon)D_02^{-m}s'$ connecting the points $\psi(\bdry\sigma)$. This gives rise to an $S\in\intcurr_1(X)$ with
$\bdry S= \Lambda_0(\bdry\Lbrack\sigma\Rbrack)$ which satisfies 
\begin{equation*}
\mass{S}\leq (1+2\varepsilon)D_02^{-m}s'\leq(1+2\varepsilon)D_0\mass{\Lbrack\sigma\Rbrack}
\end{equation*}
and
\begin{equation*}
 \spt(\eta_\#S)\subset B(\sigma, (1+2\varepsilon)^2(k+1)D_02^{-m}s').
\end{equation*}
We define $\Lambda_1\left(\Lbrack\sigma\Rbrack\right):= S$. Clearly, $\Lambda_1$ satisfy properties (i) and (ii) with $C'_1:= (1+2\varepsilon)D_0$ and
(v) with $E_1:=(1+2\varepsilon)^2(k+1)D_0$.
If $n\geq 2$ then 
\begin{equation*}
\bdry\Lambda_{n-1}\left(\bdry\Lbrack\sigma\Rbrack\right)= \Lambda_{n-2}\left(\bdry^2\Lbrack\sigma\Rbrack\right)=0
\end{equation*}
and thus the isoperimetric inequality for $\intcurr_{n-1}(X)$ and \corref{corollary:fillrad-fillvol} imply the existence of an $S\in\intcurr_n(X)$ with
$\bdry S= \Lambda_{n -1}(\bdry\Lbrack\sigma\Rbrack)$ and
\begin{equation*}
 \mass{S}\leq D_{n-1}\left[\mass{\Lambda_{n-1}\left(\bdry\Lbrack\sigma\Rbrack\right)}\right]^{\frac{n}{n-1}}
  \leq C'_n(2^{-m}s')^n\leq C'_n\mass{\Lbrack\sigma\Rbrack},
\end{equation*}
where $C'_n:= D_{n-1}\left[2n(n-1)^{\frac{n-1}{2}} C'_{n-1}\right]^{\frac{n}{n-1}}$, and
\begin{equation*}
 \spt(\eta_\#S)\subset B\left(\sigma, \left[E_{n-1}+ \mu'(1+2\varepsilon)(k+1)(n-1)^{\frac{1}{2}}(2nC'_n)^{\frac{1}{n-1}}\right]2^{-m}s'\right).
\end{equation*}

We set $\Lambda_n(\Lbrack\sigma\Rbrack):= S$ and extend $\Lambda_n$ to $G_n$ linearly and note that properties (i), (ii) and (v) are satisfied (after the 
obvious choice of $E_n$). This completes the construction of $\Lambda_n$.
In order to define $\Gamma_n$, $n\geq 1$, let again $\sigma\in\bdry\mathcal{Q}_n$ be arbitrary. Setting
\begin{equation*}
T':= \Lbrack\sigma\Rbrack-\eta_\#\Lambda_n\left(\Lbrack\sigma\Rbrack\right) - \Gamma_{n-1}\left(\bdry \Lbrack\sigma\Rbrack\right)
\end{equation*}
one easily checks that
\begin{equation*}
 \bdry T' = \bdry\Lbrack\sigma\Rbrack - \eta_\#\left(\bdry\Lambda_n\left(\Lbrack\sigma\Rbrack\right)\right)-\bdry \Gamma_{n-1}\left(\bdry \Lbrack\sigma\Rbrack\right)=0
\quad\text{and}\quad \mass{T'}\leq D''\mass{\Lbrack\sigma\Rbrack},
\end{equation*}
where $D'':=\left\{1+\left[(1+2\varepsilon)(k+1)\right]^n C'_n+2n(n-1)^{\frac{n-1}{2}} D'_{n-1}\right\}$.
Using the isoperimetric inequality for $\intcurr_n(V')$ and \corref{corollary:fillrad-fillvol} we find an $S\in\intcurr_{n+1}(V')$ with $\bdry S=T'$ and
\begin{equation*}
 \mass{S}\leq \tilde{D}_n[\mass{T'}]^{\frac{n+1}{n}}\leq 2^{-m}\tilde{D}_n(D'')^{\frac{n+1}{n}}n^{\frac{n+1}{n}}s'\mass{\Lbrack\sigma\Rbrack}
\end{equation*}
and
\begin{equation*}
 \spt(\eta_\#S)\subset B\left(\sigma, \left[E_n+ E'_{n-1}+ \tilde{\mu}'(D'')^{\frac{1}{n}}\sqrt{n}\right]2^{-m}s'\right),
\end{equation*}
where $\tilde{D}_n$ is the isoperimetric constant for $\intcurr_n(V')$ and $\tilde{\mu}'$ is the constant from \corref{corollary:fillrad-fillvol} for $V'$.
We define $\Gamma_n\left(\Lbrack\sigma\Rbrack\right):= S$ and extend it linearly to $G_n$. Clearly, $\Gamma_n$ satisfies the properties (iii), (iv) 
and (vi) with suitable choices of $D'_n$ and $E'_n$. 
This concludes the construction of the homomorphisms $\Lambda_n$ and $\Gamma_n$ for $n=0, 1,\dots, k$ with the desired properties.

Let now $\Lbrack Q\Rbrack\in\intcurr_{k+1}(V')$ be the integral current associated with $Q$ endowed with the orientation $v_1\wedge\dots\wedge v_{k+1}$ and set
$T:= \Lambda_k\left(\bdry\Lbrack Q\Rbrack\right)$. By property (ii) we have
\begin{equation*}
 \mass{T}\leq 2(k+1)C'_k k^{\frac{k}{2}}(s')^k.
\end{equation*}
Given an $S\in\intcurr_{k+1}(L^\infty(X))$ with $\bdry S= T$ we compute
\begin{align*}
 \bdry\left(\eta_\# S+\Gamma_k(\bdry\Lbrack Q\Rbrack)\right)&= \eta_\#(\bdry S) +\bdry(\Gamma_k(\bdry\Lbrack Q\Rbrack))\\
   &= \eta_\#\Lambda_k(\bdry\Lbrack Q\Rbrack)) +\bdry(\Gamma_k(\bdry\Lbrack Q\Rbrack))\\
   &= \bdry\Lbrack Q\Rbrack,
\end{align*}
where the last equality follows from property (iii).
We therefore obtain
\begin{equation*}
 \mass{\eta_\# S}\geq \mass{\Lbrack Q\Rbrack} - \mass{\Gamma_k(\bdry\Lbrack Q\Rbrack)}\geq \left[1-2^{1-m}(k+1)k^{\frac{k}{2}}D'_k\right](s')^{k+1}
\end{equation*}
and since $S$ was arbitrary we conclude
\begin{equation*}
 \fillvol_{L^\infty(X)}(T)\geq \frac{1-2^{1-m}(k+1)k^{\frac{k}{2}}D'_k}{[(1+2\varepsilon)(k+1)]^{k+1}}(s')^{k+1}.
\end{equation*}
From this and the isoperimetric inequality for $\intcurr_k(X)$ it follows that
\begin{equation*}
 \left(\frac{1-2^{1-m}(k+1)k^{\frac{k}{2}}D'_k}{D_k[(1+2\varepsilon)(k+1)]^{k+1}}\right)^{\frac{k}{k+1}}\cdot(s')^k\leq \mass{T}\leq 2(k+1)k^{\frac{k}{2}}C'_k(s')^k.
\end{equation*}
Furthermore, by the choice of $\{v_1,\dots,v_{k+1}\}$, we have
\begin{equation*}
 \fillrad(\bdry\Lbrack Q\Rbrack)\geq As'
\end{equation*}
for some constant $A>0$ only depending on $k$.
We conclude that
\begin{equation*}
 \fillrad_{L^\infty(X)}(T)\geq \frac{A-E'_k2^{-m}}{(1+2\varepsilon)(k+1)}s'.
\end{equation*}
Choose now $m\in\N$ sufficiently large to conclude the proof of the theorem.
\end{proof}

\section{Polynomial isoperimetric inequalities and the proof of \thmref{theorem:infinite-rank-isop-progression}}\label{section:proof-generalized-isop-inequaliy}

In this section we show how the decomposition theorem, \thmref{theorem:suitable-decomposition}, can be used to establish polynomial isoperimetric inequalities with exponent different from the Euclidean one. This will allow to prove \thmref{theorem:infinite-rank-isop-progression}.

The following definition generalizes the notion of cone type inequalities. 

\bd\label{definition:diameter-vol-ineq}
 Let $k\in\N$ and $\nu,\varrho>0$. A complete metric space $X$ is said to admit a diameter-volume inequality of type $(\nu,\varrho)$ for $\intcurr_k(X)$ if there
 exists $C\in(0,\infty)$ such that for every $T\in\intcurr_k(X)$ with $\bdry T=0$ and bounded support 
 \begin{equation}\label{equation:def-diam-vol-leq-1}
  \fillvol(T)\leq C\diam(\spt T)\mass{T} 
 \end{equation}
 if $\diam(\spt T)\leq 1$ and
 \begin{equation}\label{equation:nu-rho-diam-vol-ineq}
   \fillvol(T)\leq C[\diam(\spt T)]^\nu\mass{T}^\varrho 
 \end{equation}
 otherwise.
\ed

Diameter-volume inequalities of type $(1,1)$ are exactly cone type inequalities. On the other hand,
easy examples of spaces admitting diameter-volume inequalities of type $(\nu, 1)$ are simply connected homogeneous nilpotent Lie groups of class $\nu$ and, more generally, 
metric spaces all of whose subsets $B$ with $R:=\diam B<\infty$ can be contracted along curves of length at most $A\max\{R,R^\nu\}$ and which satisfy a weak 
form of the fellow traveller property (similar to that of asynchronous combings in geometric group theory). 
For precise statements in this direction see \secref{section:diameter-vol-ineq-spaces}.

\bt\label{theorem:generalized-isop-theorem}
 Let $X$ be a complete metric space, $k\in\N$, $\nu,\varrho>0$, and suppose $X$ admits a diameter-volume inequality of type $(\nu,\varrho)$ for $\intcurr_k(X)$.
 If $k=1$ set $\alpha_0:=1$. If $k\geq 2$ then suppose that $X$ admits an isoperimetric inequality of rank $\alpha_{k-1}$ 
 for $\intcurr_{k-1}(X)$ for some $\alpha_{k-1}>1$. 
 If $\nu+\varrho\alpha_{k-1}>\alpha_{k-1}$ then $X$ admits an isoperimetric inequality of rank 
 \begin{equation*}
  \alpha_k:= 1+\frac{\alpha_{k-1}}{\nu+\varrho\alpha_{k-1}-\alpha_{k-1}}
 \end{equation*}
 for $\intcurr_k(X)$ with a constant which depends only on $k, \nu,\varrho,\alpha_{k-1}$ and the constants from the isoperimetric inequality 
 for $\intcurr_{k-1}(X)$ and the diameter-volume inequality for $\intcurr_k(X)$.
\et

The definition of isoperimetric inequality of rank $\alpha$ was given in \secref{section:decomposition-theorem}. The statement of the theorem can be reformulated as follows: If $X$ admits an isoperimetric inequality for $\intcurr_{k-1}(X)$ with exponent $\mu>1$ then
\begin{equation*}
 \FV_{k+1}(r)\leq D r^{\varrho+\nu-\frac{\nu}{\mu}}
\end{equation*}
for all $r\geq 1$ and for a suitable constant $D$.
Theorems 3.4.C and 4.2.A in  \cite{Gromov-filling} and Theorem 1.2 in \cite{Wenger-GAFA} are special cases of the above theorem with
$\nu=1$ and $\alpha_{k-1}=k$. Note also, that \thmref{theorem:infinite-rank-isop-progression} stated in the introduction is a direct consequence of \thmref{theorem:generalized-isop-theorem}.


\begin{proof}[{Proof of \thmref{theorem:generalized-isop-theorem}}]
 Set $\delta:=\lambda:=1/5$ and set furthermore $\alpha:= 0$ in case $k=1$ and $\alpha:=\alpha_{k-1}$ otherwise.
 Abbreviate $F:= F_{k,\alpha,\gamma}$ and $G:= G_{k,\alpha}$, where $\gamma$ is the constant of \thmref{theorem:suitable-decomposition}.
 Let $T\in\intcurr_k(X)$ with $\bdry T=0$ and let a $R$, $T_j$ be given as in \thmref{theorem:suitable-decomposition}. Throughout this proof, the numbers
 (i) through (v) will refer to the properties listed in \thmref{theorem:suitable-decomposition}. Furthermore all the constants $E_l$  that appear
 will depend only on $k$, $\alpha$, $\gamma$ unless stated otherwise.
 Set $T_0:= R$. After possibly applying \lemref{lemma:optimal-decomposition-R} we may assume that 
 \begin{equation*}
  \diam(\spt T_j)\leq E_1 G(\mass{T_j})
 \end{equation*}
 for all $j\geq 0$ and for a constant $E_1$.
 Suppose first that $\mass{T}\leq 1$. Then, by (v), we have $\mass{T_j}\leq 3/2$ for all $j\geq 0$, so that
 \begin{equation*}
  \diam(\spt T_j)\leq E_2\mass{T_j}^{\frac{1}{k}}
 \end{equation*}
 for some constant $E_2$. It follows from the diameter-volume inequality that for each $j\geq 0$ there exists an $S_j\in\intcurr_{k+1}(X)$ 
 with $\bdry S_j=T_j$ and 
 \begin{equation*}
  \mass{S_j}\leq E_3\mass{T_j}^{\frac{k+1}{k}}
 \end{equation*}
 for some $E_3$ depending only on $k,\alpha,\gamma$ and $C_k, \nu, \varrho$. 
 This is clear if $\diam(\spt T_j)\leq 1$. If, on the other hand, 
 $\diam(\spt T_j)>1$, then we have $\mass{T_j}> E_2^{-k}$ and from this the inequality readily follows. Finally, we have
 \begin{equation*}
  \sum_{j=0}^\infty\mass{S_j}\leq E_3\sum_{j=0}^\infty\mass{T_j}^{\frac{k+1}{k}}\leq 
   E_3\left[\sum_{j=0}^\infty\mass{T_j}\right]^{\frac{k+1}{k}}\leq E_3\left[\frac{3}{2}\right]^{\frac{k+1}{k}}
      \mass{T}^{\frac{k+1}{k}}.
 \end{equation*}
 Therefore, $\sum_{j=0}^nS_j$ is a Cauchy sequence with respect to mass and therefore converges to some $S\in\intcurr_{k+1}(X)$,
 which clearly satisfies $\bdry S=T$ and 
 \begin{equation*}
  \mass{S}\leq \sum_{j=0}^\infty\mass{S_j}\leq E_3\left[\frac{3}{2}\right]^{\frac{k+1}{k}}\mass{T}^{\frac{k+1}{k}}.
 \end{equation*}
 This proves the theorem if $\mass{T}\leq 1$. In case $\mass{T}>1$ define
 $J:=\{j\geq 0: \diam(\spt T_j)>1\}$. If $j\in J$ then, by the diameter-volume inequality, there exists an $S_j\in\intcurr_{k+1}(X)$ with 
 $\bdry S_j=T_j$ and 
 \begin{equation*}
  \mass{S_j}\leq C_k \mass{T_j}^{\frac{\nu}{\alpha}+\varrho}.
 \end{equation*}
 On the other hand, if $j\not\in J$ then, again by the diameter-volume inequality, there exists an $S_j\in\intcurr_{k+1}(X)$ with $\bdry S_j=T_j$ and
 \begin{equation*}
  \mass{S_j}\leq C_k \diam(\spt T_j) \mass{T_j}\leq C_k\mass{T_j}.
 \end{equation*}
 Since $\nu+\varrho\alpha>\alpha$ and $\mass{T}>1$ we have
 \begin{equation*}
  \begin{split}
  \sum_{j=0}^\infty\mass{S_j} &\leq C_k\sum_{j\in J}\mass{T_j}^{\frac{\nu}{\alpha}+\varrho} + C_k\sum_{j\not\in J}\mass{T_j}\\
           &\leq C_k\left[\sum_{j\in J}\mass{T_j}\right]^{\frac{\nu}{\alpha}+\varrho} + \frac{3}{2}C_k\mass{T}\\
           &\leq 2C_k\left[\frac{3}{2}\right]^{\frac{\nu}{\alpha}+\varrho}\mass{T}^{\frac{\nu}{\alpha}+\varrho}.
  \end{split}
 \end{equation*}
 It now follows exactly as above that $\sum_{j=0}^nS_j$ converges in mass to some $S\in\intcurr_{k+1}(X)$ which has the desired properties.
 This concludes the proof.
\end{proof}

\section{Appendix: Metric spaces admitting cone type and diameter-volume inequalities}\label{section:diameter-vol-ineq-spaces}

Let $(X,d)$ be a metric space, $B\subset X$ with $\diam B<\infty$, and $h,H>0$.
Suppose there exists a Lipschitz map $\varphi: [0,1]\times B\to X$ with the following properties:
\begin{enumerate}
 \item There exists $x_0\in X$ such that $\varphi(0,x)= x_0$ and $\varphi(1,x)=x$ for all $x\in B$. 
 \item The lengths of the curves $t\mapsto \varphi(t,x)$, $x\in B$, are bounded above by $h$.
 \item For every $x\in B$ there exists a relatively open neighborhood $U_x\subset B$ of $x$ and a continuous family $\varrho_{x'}$, $x'\in U$, 
 of reparametrizations of $[0,1]$ such that
   $\varrho_x(t)=t$ and
   \begin{equation*}
    d(\varphi(t,x), \varphi(\varrho_{x'}(t), x'))\leq H d(x,x')\quad\text{ for all $x'\in U_x$ and all $t\in[0,1]$.}
   \end{equation*}
\end{enumerate}
 Here, a map $\nu:[0,1]\to [0,1]$ is said to be a reparametrization of $[0,1]$ if it is continuous, non-decreasing and satisfies $\nu(0)=0$ and $\nu(1)=1$.
 We call $\varphi$ as above a Lipschitz contraction of $B$ with parameters $(h,H)$.
 
 \bd
  Let $X$ be a metric space and  $h,H:[0,\infty)\to[0,\infty)$ continuous functions. If
every subset $B\subset X$ with 
 $R:=\diam B<\infty$ has a Lipschitz contraction with parameters $(h(R), H(R))$ then $X$ is said to admit generalized combings with length function $h$ and distortion function
 $H$.
\ed

Simple examples of spaces admitting generalized combings are given as follows: Normed spaces, ${\rm CAT}(0)$-spaces and, more generally, geodesic metric spaces with convex metric (see \secref{section:asymptotic-rank} for the definition) admit generalized combings with length function $h(R):=  R$ and distortion function $H(R)=1$. Simply connected homogeneous nilpotent Lie groups of class $c$, endowed with a left-invariant Riemannian metric, admit generalized combings with length function $$h(R):= A\max\{R, R^c\}$$ and
   distortion function $H(R)\equiv L$ for constants $A,L$, see \cite{Pittet-homogeneous-nilpotent}.
We now prove the following proposition.
 \bp\label{proposition:coneineq-locally-contractible}
  Let $X$ be a complete metric space, $T\in \intcurr_k(X)$ a cycle with bounded support and $h,H>0$. Suppose there exists a Lipschitz contraction $\varphi$ of $B:=\spt T$ 
  with parameters $(h,H)$. Then there exists $S\in\intcurr_{k+1}(X)$ with $\bdry S=T$ and such that
  \begin{equation*}
   \mass{S}\leq [k(k+1)]^{\frac{k}{2}}hH^k\mass{T}.
  \end{equation*}
  \ep
Before turning to the proof we mention the following immediate consequence of the proposition. 
Given $\nu\geq 1$, $\mu\geq 0$ and $A, L>0$ define
\begin{equation*}
 h_{\nu}(r):=\left\{\begin{array}{l@{\qquad}l}
  Ar & 0\leq r\leq 1\\
  Ar^\nu & 1 < r <\infty
 \end{array}\right. \qquad\text{and}\qquad H_{\mu}(r):=\left\{\begin{array}{l@{\qquad}l}
   L & 0\leq r\leq 1\\
   Lr^{\,\mu} & 1 < r <\infty.
 \end{array}\right.
\end{equation*}

We have:

\bc
 Let $X$ be a complete metric space. If $X$ admits generalized combings with length function $h_{\nu}$ and distortion function $H_{\mu}$ then
 $X$ admits a diameter-volume inequality of type $(\nu+k\mu, 1)$ for $\intcurr_k(X)$ for every $k\geq 1$. In particular, if $\nu=1$ and $\mu=0$ then $X$ admits a cone type inequality for $\intcurr_k(X)$.
\ec

For the definition of volume-diameter inequalities see \secref{section:proof-generalized-isop-inequaliy}.

 \begin{proof}[{Proof of \propref{proposition:coneineq-locally-contractible}}]
  By Theorem 4.5 of \cite{Ambr-Kirch-curr} it is enough to consider the case $T=\psi_{\#}\Lbrack\theta\Rbrack$ for a biLipschitz map $\psi: K\subset\R^k\to X$
  and $\theta\in L^1(K,\Z)$. Set $B:=\spt T$ and let $\varphi$ be a Lipschitz contraction of $B$. Set $S:= \varphi_{\#}([0,1]\times T)$, where the product $[0,1]\times T$ of an interval with $T$ is defined as in \cite{Wenger-GAFA}, and note that $S\in\intcurr_{k+1}(X)$ and $\bdry S=T$.
  Let $f$ and $\pi_1,\dots,\pi_{k+1}$ be 
  Lipschitz functions on $X$ with $f$ bounded and such that $\lip(\pi_i)\leq 1$ for all $i$.
  We define $\tilde{\varphi}(t,z):= \varphi(t,\psi(z))$ for $t\in [0,1]$ and $z\in K$ and $\pi:= (\pi_1,\dots, \pi_{k+1})$.
  Let $(t,z)\in[0,1]\times K$ be such that $\pi\circ\tilde{\varphi}$ is differentiable at $(t,z)$ with non-degenerate differential, which we denote by $Q$, 
  $\psi$ is metrically differentiable at $z$ in the sense of \cite{Kirchheim} and the curve 
  $\gamma_z(t):=\tilde{\varphi}(t,z)$ is metrically differentiable at $t$. We may assume without loss of generality that $\pi\circ\tilde{\varphi}(t,z)=0$.
  Denote by $P$ the orthogonal projection of $\R^{k+1}$ onto the orthogonal complement of $Q(\R\times\{0\})$.
  We claim that 
  \begin{equation*}
   \|P(Q(0,v))\|\leq \sqrt{k+1}H\md\psi_z(v)
  \end{equation*}
  for all $v\in\R^k$.
  In order to see this fix $v\in\R^k\ohne\{0\}$ and choose, for each $r>0$ sufficiently small, some $t_r$ with $t=\varrho_{\psi(z+rv)}(t_r)$, where $\varrho_{x'}$ denotes the 
  family of reparametrizations of $[0,1]$ around $\psi(z)$. It is easy to see that $|t_r-t|\leq Cr$ for some constant $C$ and all $r>0$ sufficiently small. It then follows that
  \begin{equation*}
   \begin{split}
    \|P(Q(0,v))\| &= \lim_{r\searrow 0}\frac{1}{r}\|P(\pi\circ\varphi(t,\psi(z+rv)) - \pi\circ\varphi(t_r, \psi(z))) + P(\pi\circ\varphi(t_r, \psi(z)))\|\\
      & \leq \lip(\pi)\limsup_{r\searrow0}\frac{1}{r} d(\varphi(t, \psi(z+rv)), \varphi(t_r, \psi(z)))\\
      &\leq \sqrt{k+1}H\md\psi_z(v).
   \end{split}
  \end{equation*}
 This proves the claim and furthermore yields
 \begin{equation*}
  \left|\det Q\right| \leq (k+1)^\frac{k}{2} H^k\jac_1(\md(\gamma_z)_t)\jac_k(\md\psi_z).
 \end{equation*}
 We use this, the area formula in \cite{Kirchheim} and Lemma 9.2 and Theorem 9.5 in \cite{Ambr-Kirch-curr} to conclude
  \begin{equation*}
   \begin{split}
    |S&(f,\pi_1,\dots,\pi_{k+1})|\\ &\leq\int_{[0,1]\times K} \left|\theta(z) f(\tilde{\varphi}(t,z)) \det(D_{(t,z)}(\pi\circ\tilde{\varphi}))\right| 
     d\lm^{k+1}(t,z)\\
       &\leq (k+1)^{\frac{k}{2}}H^k\int_{[0,1]\times K}|\theta(z)f(\varphi(t,\psi(z)))|\jac_1(\md(\gamma_z)_t)\jac_k(\md\psi_z)d\lm^{k+1}(t,z)\\
       &\leq [k(k+1)]^{\frac{k}{2}}H^k\int_{[0,1]\times X}|f(\varphi(t,x))|\jac_1(\md(\gamma_{\psi^{-1}(x)})_t)d(\lm^1\times\|T\|)(t,x)\\
   \end{split}
  \end{equation*}
  and thus $$\|S\|\leq [k(k+1)]^{\frac{k}{2}}H^k\varphi_{\#}[g(t,x)d(\lm^1\times\|T\|)]$$ with $g(t,x):= \jac_1(\md(\gamma_{\psi^{-1}(x)})_t)$. This 
  completes the proof.
  \end{proof} 
We finally mention that diameter-volume inequalities can also be established for spaces with nice local geometry on which asynchronously combable groups with polynomial 
length functions act properly and cocompactly by isometries. See for example Chapter 10 of \cite{Epstein-et-all} or Part II of \cite{Riley-filling-notes}.

\end{document}